\documentclass{amsart}

\usepackage{amscd}
\usepackage{amssymb}
\usepackage[arrow,matrix,curve,ps]{xy}
\xyoption{dvips}
\CompileMatrices

\newcommand{\CZ}{\mbox{$\mathbb{C}$}}
\newcommand{\QZ}{\mbox{$\mathbb{Q}$}}
\newcommand{\ZZ}{\mbox{$\mathbb{Z}$}}
\newcommand{\NZ}{\mbox{$\mathbb{N}$}}
\newcommand{\Pone}{\mbox{$\mathbb{P}_1$}}
\newcommand{\Ptwo}{\mbox{$\mathbb{P}_2$}}
\newcommand{\Pthree}{\mbox{$\mathbb{P}_3$}}
\newcommand{\Pfive}{\mbox{$\mathbb{P}_5$}}
\newcommand{\name}[1]{{\sc #1}\index{#1}}

\DeclareMathOperator{\Aut}{Aut}
\DeclareMathOperator{\Pic}{Pic}
\DeclareMathOperator{\Sing}{Sing}

\theoremstyle{plain}
\begingroup
\newtheorem{thm}{Theorem}[section]

\newtheorem{lem}[thm]{Lemma}
\newtheorem{prop}[thm]{Proposition}
\endgroup

\theoremstyle{remark}
\newtheorem{assumption}[thm]{Assumption}
\newtheorem{example}[thm]{Example}
\newtheorem{notation}[thm]{Notation}
\numberwithin{equation}{section}

\author{Stefan Kebekus}
\date{\today}
\email{stefan.kebekus@uni-bayreuth.de}
\address{Stefan Kebekus\\ Mathematisches Institut der Universit\"at
Bayreuth\\ 95440 Bayreuth\\ Germany\\ FAX: +49 (0)921/55-2785 }

\date{April 20, 1998}

\title{Relatively Minimal Quasihomogeneous Projective 3-Folds}
\thanks{The author was supported by scholarships of the
Graduiertenkollegs ``Geometrie und mathematische Physik'' and
``Komplexe Mannigfaltigkeiten'' of the Deutsche
Forschungsgemeinschaft}

\begin{document}

\begin{abstract}
In the present work we classify the relatively minimal 3-dimensional
quasihomogeneous complex projective varieties under the assumption
that the automorphism group is not solvable. By relatively minimal we
understand varieties $X$ having at most $\QZ$-factorial terminal
singularities and allowing an extremal contraction $X\to Y$ where
$\dim Y <3$.

\end{abstract}

\maketitle
\tableofcontents

\section{Introduction}

Let $X$ be a smooth projective threefold and $G$ a connected algebraic
group acting algebraically on $X$. In this context all steps of the
minimal model program are equivariant. If one assumes additionally
that $G$ acts almost transitively, which is to say that the $G$-action
has an open orbit, then the minimal model program always leads to a
contraction of fiber type over a base $Y$, i.e.~to a relatively
minimal model (see e.g.~\cite{K98a} for details on this). The
relatively minimal models of smooth varieties are always
$\QZ$-factorial. Here we classify these varieties under the assumption
that $G$ is linear algebraic and not solvable.

We now list the non-trivial examples which occur in the
classification. Notation: a ``linear bundle'' is a variety of the form
$\mathbb{P}(E)$, where $E$ is a vector bundle. Call $\mathbb{P}(E)$
``splitting'' if $E$ is.

\begin{itemize}
\item The special \name{Fano} varieties $V_5$ and $V^S_{22}$ described by 
\name{Mukai} and \name{Umemura} ---see \cite{MU83}.

\item The weighted projective spaces $\mathbb P_{(1,1,2,3)}$ and
$\mathbb P_{(1,1,1,2)}$. The first space is described in detail in
\cite[ex.~4.1]{K98c}, the latter is blow-down of the negative section
of $\mathbb P(\mathcal O_{\Ptwo}(2)\oplus\mathcal O_{\Ptwo})$.

\item Varieties over $Y\cong\Pone$ which are locally isomorphic to a
deformation of a quadric surface, and certain quotients of these
varieties by $\ZZ_2$. They are described in detail in the sections
\ref{V3ex} and \ref{V4ex} of the present paper and called the
``quadric- and $\mathbb{F}_4$-degenerations''.

\item Singular varieties arising as quotients by $\ZZ_2$ of a splitting 
linear $\Pone$-bundle over $Y\cong\Ptwo$; see
example~\ref{Surf:SingSurf}. Abusing language, call these the
``singular $\Pone$-bundles over $Y\cong \mathbb F_2$''.

\item Linear $\Pone$-bundles over \name{Hirzebruch}-surfaces $Y$ which are 
constructed in sections~\ref{ex01}--\ref{Surf:X_Sigma_0_n} by starting
with a trivial $\Pone$-bundle and repeatedly performing certain
elementary transformations; name these varieties the ``diagonally
twisted bundles''.
\end{itemize}

\begin{thm}\label{mainthm}
Let $X$ be a 3-dimensional projective complex variety with at most
$\QZ$-factorial terminal singularities and $G$ be a connected linear
algebraic group acting algebraically and almost transitively on $X$ so
that the kernel of $G\to \Aut(X)$ is discrete. Assume that $G$ is not
solvable. If there exists an extremal contraction $\phi:X\rightarrow
Y$ with $\dim Y<3$, then $X$ is isomorphic to one of the following:
\begin{itemize}

\item If $Y$ is a point, then $X\cong \mathbb P_{(1,1,2,3)}$, 
$\mathbb P_{(1,1,1,2)}$, or $X$ is one of the following \name{Fano}
varieties: $\Pthree$, $\QZ_3$, the 3-dimensional quadric, $V_5$ or
$V_{22}^S$.

\item If $\dim Y=1$, $X$ is one of the quadric- and
$\mathbb{F}_4$-degenerations, or a linear $\Ptwo$-bundle over $\Pone$.

\item If $\dim Y=2$ and if $Y$ is singular, $X$ is a singular 
$\Pone$-bundle over $\mathbb F_2$. Otherwise, $X\cong Y\times\Pone$,
where $Y$ is an arbitrary $G$-quasihomogeneous surface, or $X$ is
smooth and one of the following holds:
\begin{itemize}
\item  $X$ is one of the diagonally twisted bundles, or a splitting 
linear bundle and $Y$ is a \name{Hirzebruch}-surface $\Sigma_n$.

\item $X$ is the full flag variety $F_{(1,2)}(3)$, a splitting linear 
bundle, a quotient of one of the models over $\Sigma_0$ or a blow-down 
of the diagonally twisted bundle $X_{\Sigma_1,k_0,0}$. In all these
cases $Y\cong \Ptwo$.
\end{itemize}
\end{itemize}
\end{thm}

We underline that theorem~\ref{mainthm} is a classification of the
{\em relatively minimal} models.

Although we found it easier to use the dimension of $Y$ to structure
the present paper, it might be worth while to briefly discuss the
classification based on the dimension of the generic orbits of a
maximal semi-simple subgroup $S$ of $G$.

If $S$ acts almost transitively, the case of primary interest is that
where $S\cong SL_2$.  Here $X$ must in fact be smooth (see
lemma~\ref{factlem}). In this setting the case that $\dim Y=0$ has
been treated in the literature (\cite{MU83}, and the papers of
\name{Iskovskih}). In the other cases where $\dim Y=1$ or $2$, one 
could apply the methods and results of \cite{MJ90} if one would extend
this to all possible isotropy groups. We choose a different approach
and construct all varieties explicitly.

If the generic $S$-orbit is 2-dimensional, the models over surfaces can
be easily described. If $Y$ is a curve, we again give an explicit
construction of all the possible varieties ---locally these are the
well-known deformations of the cones over rational normal curves of
degree 2 or 4. The remaining case where $\dim Y=0$ is slightly more
involved and requires a line of argumentation that does not fit well
into the present work. Thus, we have chosen to treat this case in a
different paper \cite{K98c}.

Finally, if the generic $S$-orbit is $1$-dimensional, then $X$ is a
product $Y\times \Pone$.

The author would like to thank A.~Huckleberry, and T.~Peternell for
support and valuable discussions. Part of the work on this paper was
carried out during our visit to the University of Grenoble. We are
thankful to the members of that institute for their kind hospitality.

\section{Minimal Models over curves}\label{MinOverCurve} 

In this section we consider the following situation:
\begin{assumption}\label{curve:assumption}
Let $X$ and $G$ be as in theorem~\ref{mainthm} and $\phi:X\rightarrow
Y$ be an extremal contraction to a curve.
\end{assumption}

Recall that $Y$ is necessarily normal and quasihomogeneous with
respect to an algebraic action of the linear algebraic group
$G$. Thus, $Y\cong\Pone$. If $X_\eta$ is a general $\phi$-fiber, it is
\name{del Pezzo} and quasihomogeneous. Therefore it is isomorphic to
either
\begin{itemize}
\item the projective plane $\mathbb{P}_2$
\item $\Sigma_0 \cong \mathbb{P}_1 \times \mathbb{P}_1$
\item the first \name{Hirzebruch}-surface $\Sigma_1$, or
\item a blow-up of $\Sigma_0$ in at most two points $x_1$ and $x_2$
such that both natural projections $\pi_i:\Sigma_0 \rightarrow \Pone$
satisfy $\pi_i(x_1)\not = \pi_i(x_2)$.
\end{itemize}
We will show that only the first two cases occur. To start with, fix
some notation:
\begin{notation}
Under the above assumptions, for $\eta \in Y$ let
$X_\eta=\phi^{-1}(\eta)$ be the associated fiber and $G_\eta$ be the
stabilizer of $X_\eta$, i.e.~the isotropy group of $\eta$.
\end{notation}

The following simple observation is crucial and will be constantly
used in the sequel:
\begin{lem}[Homology Lemma]\label{homology_lemma}
Let $\phi:X\rightarrow Y$ be an extremal contraction ($Y$ not
necessarily a curve), $D\in Div(X)$ an irreducible divisor and $y\in
Y$ a point. If $D\cap X_\eta$ is a nontrivial effective divisor, then
it intersects every curve in $X_\eta$ positively.
\end{lem}
\begin{proof}
There is a curve $C\subset X_\eta$ intersecting $D$ in a finite
set. So $C.D>0$. Let $C'\subset X_\eta$ be any other curve. Since
$\phi$ is a contraction, there exist $a,b \in \QZ^+$ such that
$a[C]=b[C']$ as homology classes. Thus $D.C'=\frac{a}{b}D.C>0$.
\end{proof}

Now we can characterize the $\phi$-fibers:
\begin{lem}
Under assumptions~\ref{curve:assumption}, a generic fiber $X_\eta$ is
isomorphic to $\Ptwo$ or to a 2-dimensional quadric.
\end{lem}
\begin{proof}
Suppose to the contrary. Then there exist (-1)-curves in
$X_\eta$. Choosing one of them, say $C$, then $D:=\overline{G.C}$ is a
divisor intersecting $X_\eta$ in $G_\eta.C$, i.e.~a finite number of
(-1)-curves. We will treat the possibilities for $X_\eta$ separately
and show that in each case the existence of $D$ yields a contradiction
to the homology lemma.

Assuming $X_\eta \cong \Sigma_1$, there is a unique (-1)-curve $C$.
Hence $D\cap X_\eta = C$, and there are curves $C'\subset X_\eta$ with
$C'.D =0$, a contradiction!

If $X_\eta \cong \Pone\times\Pone$ blown up in one point, there are
three (-1)-curves $C_1$, $C_2$ and $C_3$ contained in $X_\eta$. They
satisfy $C_1.C_2 = C_2.C_3 = 1$ and $C_1.C_3 = 0$. Set
$D:=\overline{G.C_2}$. Then, if $D\cap X_\eta$ contains $C_1$ and
$C_2$, $C_1.D =0$. If $D\cap X_\eta$ contains $C_2$ and $C_3$, $C_3.D
=0$. As a last possibility, $D\cap X_\eta = C_2$. Then there exists a
curve in $X_\eta$ which does not intersect $D$ at all. In any case,
the homology lemma is violated.

The last case is that $X_\eta \cong \Pone\times\Pone$ blown up in two
points as described above. First, we remark that a 1-dimensional
subgroup $H<G$ acting non-trivially on $Y$ cannot be isomorphic to
$\CZ$: if it were, since it's isotropy at a generic point $\eta\in Y$
is trivial, given any (-1)-curve $C\subset X_\eta$,
$D:=\overline{H.C}$ would be a divisor, $D\cap X_\eta = C$, and there
would exist curves in $X_\eta$ not intersecting $D$. In particular,
this implies that $G$ acts as $\CZ^*$ on $Y$.

On the other hand, since $\Aut^0(X_\eta)\cong \CZ^*\times \CZ^*$, it
follows that $G$ acts as a torus $(\CZ^*)^3$! A contradiction to the
assumptions.
\end{proof}

\subsection{The Construction of the Quadric-Degenerations}\label{V3ex} 

In this section we construct models over $\Pone$ whose generic fibers
are smooth quadrics. For this, consider the space $V_0 := \Pthree
\times \CZ$ with coordinates $([x:y:z:w],\lambda)$. For any odd
integer $k>0$, let $X_0^k$ be the quasi-projective variety given by:
$$
X_0^k := \{([x:y:z:w],\lambda)\in \Pthree\times \CZ:
4xz-y^2=\lambda^kw^2 \}
$$
Let $SL_2$ act on $\Pthree$ via a direct sum of the one- and
three-dimensional irreducible representations, thus stabilizing the
quadric $\{4xz-y^2=\lambda^kw^2\}$. The group $H^*\cong \CZ^*$ acts as
follows: 
$$
\xi ([x:y:z:w],\lambda) = ([x:y:z:\xi^{-k}w],\xi^2 \lambda).
$$
A direct calculation shows that $G:=H^*\times SL_2$ acts and
stabilizes $X_0^k$.

Choosing another odd number $l$, we construct a similar
quasi-projective variety $X^l_\infty$ over $\CZ$: Again $V_\infty :=
\Pthree\times\CZ$ and  $X^l_\infty := \{ 4xz-y^2 = \lambda^lw^{2}
\}$. Let $SL_2$ act as above and let $H^*$ act by:
$$
\xi: ([x:y:z:w],\lambda) \mapsto([x:y:z:\xi^l w],\xi^{-2} \lambda).
$$

The last step of the construction consists in gluing $V_0$ and
$V_\infty$ in order to obtain a $\Pthree$-bundle over $\Pone$
which contains the desired quasihomogeneous space. Define the
equivalence relation
\begin{eqnarray*}
V_0 \ni ([x_0:y_0:z_0:w_0],\lambda_0) \sim
([x_\infty:y_\infty:z_\infty:w_\infty],\lambda_\infty) \in
V_\infty  \\ 
: \Leftrightarrow \lambda_0 \lambda_\infty = 1 \text{\ and \ }
[x_0:y_0:z_0:w_0]=[x_\infty:y_\infty:z_\infty:w_\infty
\lambda^{(k+l)/2}_\infty].
\end{eqnarray*}
Consider the equation defining $X^k_0$ and substitute the equivalent
coordinates of $V_\infty$:
\begin{eqnarray*}
4x_0z_0-y_0^2 &=& \lambda^k_0 w_0^2 \\
\Leftrightarrow 4x_\infty z_\infty-y_\infty^2 &=&
\frac{1}{\lambda^k_\infty} (w_\infty \lambda^{(k+l)/2}_\infty)^2 \\ 
\Leftrightarrow 4x_\infty z_\infty^2-y_\infty^2 &=&
w_\infty^2\lambda^l_\infty 
\end{eqnarray*}
the last equation is that which defines $X^l_\infty$. 

There are several things to show:

\subsubsection{$X^{(k,l)}$ has $\QZ$-factorial terminal singularities}
As a first step, we claim that $H^*$ acts trivially on the divisor
class group $Cl(X^k_0)$. For this, note that if $\tilde X^{(k,l)}$ is
an $H^*$-equivariant resolution of the singularities, then there
exists an $H^*$-equivariant surjection $Cl(X^{(k,l)})\rightarrow
Cl(X^k_0)$ and an $H^*$-equivariant injection
$Cl(X^{(k,l)})\rightarrow Cl(\tilde X^{(k,l)})=\Pic(\tilde
X^{(k,l)})$. But every component of $\Pic(\tilde X^{(k,l)})$ is a
compact torus, i.e.~the only algebraic $H^*$-action is trivial ---see
\cite[Lect.~19ff]{Mum66} for the fact that the action of $H^*$ on
$\Pic(\tilde X^{(k,l)})$ is algebraic.

Second, observe that $X^k_0$ has an isolated cDV singularity at
$([0:0:0:1],0)$, which is terminal of index one
(cf.~\cite[par.~1]{Reid83}). Furthermore, $X^1_0$ is smooth!  

We claim that all divisors $D\subset X_k$ are
$\QZ$-\name{Cartier}. Define a map $\gamma: X^k_0 \rightarrow X^1_0$
by $\gamma: ([x,y,z,w],\lambda)\mapsto ([x,y,z,w],\lambda^k)$. This is
a quotient of $X^k_0$ by an action of $\ZZ_k$. Observe that $$
D':=\sum_{\xi\in H^*, \xi^{2k}=1}\xi D $$ is $\ZZ_k$-invariant, hence
\name{Cartier}. This is the place where we need ``$k$ odd''. As $H^*$
acts trivially on the divisor class group of $X^k_0$, $D'$ is linearly
equivalent to a multiple of $D$. Consequently, $D$ is
$\QZ$-\name{Cartier} indeed.

The same argumentation holds for $X^l_\infty$.

\subsubsection{$X^{(k,l)}$ is $G$-quasihomogeneous}
In order to see that the group actions on the quasi-projective pieces
extend to the entire variety, we show that if $v_0 =
([x_0:y_0:z_0:w_0],\lambda_0) \sim
([x_\infty:y_\infty:z_\infty:w_\infty],\lambda_\infty) = v_\infty$ and
$g \in G$, then $g.v_0 \sim g.v_\infty$.

A simple calculation shows that this holds if $g\in SL_2$. Similarly,
if $\xi\in H^*$,
\begin{eqnarray*}
\xi([x_0:y_0:z_0:w_0],\lambda_0) &=& ([x_0: y_0: z_0 : \xi^{-k} w_0],
\xi^2 \lambda_0)\\ 
\xi([x_\infty:y_\infty:z_\infty:w_\infty],\lambda_\infty)
&=&([x_\infty:y_\infty:z_\infty:\xi^l
w_\infty],\xi^{-2} \lambda_\infty) 
\end{eqnarray*}
Now note that $(\lambda_0\xi^2)(\lambda_\infty\xi^{-2}) =
\lambda_0\lambda_\infty$ and $\xi^l w_\infty(\lambda_\infty
\xi^{-2})^{(k+l)/2} = w_0\xi^{-k}$, showing that
$\xi([x_0:y_0:z_0:w_0],\lambda_0) \sim
\xi([x_\infty:y_\infty:z_\infty:w_\infty],\lambda_\infty)$. Due to the
product structure, $g.v_0 \sim g.v_\infty$ for all $g \in G$.
 
\subsubsection{$X^{(k,l)}$ can be \name{Mori}-contracted to
$\Pone$} Perform a relative \name{Mori} contraction
$\psi:X\rightarrow Z$ over $\Pone$. Note that if $X_\mu$ is an
arbitrary fiber of the map $X\rightarrow \Pone$, then all curves
contained in $X_\mu$ are equivalent as homology cycles: this is clear
for the singular fibers over $0$ and $\infty$ because they are
singular quadrics, and also true for the generic fibers because the
action of $\pm 1\in H^*$ swaps horizontal and vertical
directions. Consequently, $\psi(X)=\Pone$, and the claim is shown.

\subsection{The Characterization of Quadric-Degenerations}
We will show that every model over $\Pone$ whose generic fiber
$X_\eta$ is a quadric is isomorphic to some $X^{(k,l)}$. Identify
$X_\eta$ with $\mathbb{P}_1 \times \mathbb{P}_1$, and let $\pi_1$ and
$\pi_2:X_\eta \rightarrow \mathbb{P}_1$ be the standard
projections. Call $\pi_1$-fibers ``vertical'' and $\pi_2$-fibers
``horizontal''.

\begin{prop}[Characterization of Quadric-Degenerations]\label{PropV3ex} 
Under assumptions~\ref{curve:assumption}, if the generic fiber
$X_\eta$ is isomorphic to a 2-dimensional quadric, then $X$ is
isomorphic to one of the quadric-degenerations constructed in
section~\ref{V3ex}.
\end{prop}
We subdivide the proof into a number of steps:
\begin{proof}[Step 1: Description of the $S$-Action] 
Let $C\subset X_\eta$ be a horizontal curve and $H<G$ a one-parameter
group acting non-trivially on $Y$. Let $H_\eta<H$ be the stabilizer of
$X_\eta$ and set $E:=\overline{H.C}$. If $E_\eta := E\cap X_\eta =
H_\eta.C$ is a union of finitely many horizontal curves, then $E_\eta$
does not intersect a general horizontal curve, contradicting the
homology lemma~\ref{homology_lemma}. Thus there exists $h\in H_\eta:
h.C$ is not horizontal. In particular, $H_\eta$ is not trivial and
$H\cong\CZ^*$ is a torus.

Let $S$ be a maximal semi-simple subgroup of $G$. As all one-parameter
subgroups of $G$ acting non-trivially on the base are necessarily
tori, $S$ acts trivially on $Y$. If some $S' \cong SL_2$ in $S$ would
act only on one factor of $X_\eta\cong \Pone\times\Pone$, we derive a
contradiction as follows: let $T<S'$ be a maximal torus and $F\subset
X$ it's fixed point set. Since $X_\eta$ was chosen to be a general
$\phi$-fiber, the $S'$-orbits in the neighboring fibers are
1-dimensional, too. So $F$ is a divisor. By assumption, $F\cap X_\eta$
is the union of two horizontal (or vertical) curves, a contradiction
to the homology lemma. Thus $S=SL_2$ and it's action is diagonal. In
particular, there exists an $S$-invariant diagonal $\delta\subset
X_\eta$.
\end{proof}
\begin{proof}[Step 2: The Embedding into a Linear Bundle]
We claim $\delta$ is also invariant under $G_\eta$. Assume to the
contrary and let $g\in G_\eta$ be an element not stabilizing
$\delta$. But $S$ has only two orbits in $X_\eta$, namely $\delta$ and
$X_\eta\setminus \delta$, so that any group containing $S$ and $g$
acts transitively on $X_\eta$, i.e.~contains $SL_2\times SL_2$. But
this is absurd, as we have seen.

Set $D:=\overline{H.\delta}$. The desingularization $\tilde{D}$ of $D$
is then quasihomogeneous. By classification,
$\tilde{D}\cong\Pone\times\Pone$. This has two consequences: First, as
$S$ acts transitively on the fibers of $\tilde{D} \rightarrow Y$, all
fibers of the map $\tilde{D} \rightarrow D$ are discrete, and the
$S$-action on $D$ does not have a fixed point.  Since the
singularities of $X$ are isolated, $D$ does not meet the singular set
of $X$. Thus, $D$ is \name{Cartier}. Second, if $U<S$ is a unipotent
group, then the $U$-fixed points in $D$ form a curve which is 
mapped injectively onto $\Pone$. This already shows that $\phi$ has
maximal rank along these curves so that there are no multiple
fibers. Furthermore, is $X_\nu$ is {\em any} $\phi$-fiber, then
$X_\nu$ is smooth along $D \cap X_\nu$. 

Recall that $\phi$ being an extremal contraction implies that $D$ is
relatively ample. As $X_\nu \cap D$ is ample and $S$-invariant, there
is no $S$-invariant curve in $X_\nu \cap D$. In particular, the
singular set of $X_\nu$ is discrete. Since all fibers are
\name{Cohen}-\name{Macaulay}, it follows from \name{Serre}'s criterion 
that they are normal. But the only normal $SL_2$-surfaces containing
an ample $S$-invariant divisor of self-intersection 2 are the
2-dimensional quadrics. Thus, we conclude that the dimension of the
linear system $|X_\nu \cap D|$ is independent of $\nu\in Y$ and that
$D$ is relatively very ample, i.e.~there exists an embedding
$X\to\mathbb P(E)$, where $E:=\phi_*(\mathcal{O}(D))$ is a rank
4-vector bundle.
\end{proof}
\begin{proof}[Step 3: Local Description]
Knowing that the intersection of $D\cap X_\eta$ yields
an equivariant embedding $X_\eta\rightarrow \Pthree$, one sees that
there is an $S$-stable splitting $E=E_3\oplus E_1$, where $E_3$ is of
rank three and $S$ acts on the fibers via it's irreducible
3-dimensional representation and $E_1$ is 1-dimensional with trivial
$S$-action. Let $T<SL_2$ be the diagonal matrices. Then the direct sum
decomposition of the irreducible $SL_2$-representations into
$T$-weight spaces yields a $T$-stable splitting $E_3 = E^{-2}_3 \oplus
E^{0}_3 \oplus E^{2}_3$, where $T$ acts on the total space of $E^i_3$
with weight $i$.

As a next step, choose a $G$-invariant affine subset $\CZ\cong U^0
\subset Y\cong \Pone$ containing one of the $G$-fixed points in
$Y$. Let $y$ be a bundle coordinate for $E^0_3$ over $U^0$; we view
that as giving a $T$-equivariant map from  $E^0_3$ into the
standard 3-dimensional $SL_2$-representation space $V_2$. In order to
obtain an $SL_2$-equivariant map $E_3|_{U^0}\rightarrow V_2$,
conjugate $y$ with the going-up and going-down operators in
$SL_2$. This way we obtain coordinates $x$ and $z$ for $E^{-2}_3$ and
$E^2_3$, respectively, giving the desired map to $V_2$.

Use these coordinates to view $X^0:=\phi^{-1}(U^0)$ as a subset of
$\Pthree\times \CZ$. The generic fiber is an $S$-invariant quadric,
hence given by $c(4xz-y^2)=c'w^2$ where $c$, $c'\in \CZ^*$. Thus,
after appropriate choice of coordinates, $X\cap\phi^{-1}(U_0)$ is
given by $4xz-y^2=\lambda^kw^2$ or $\lambda^k(4xz-y^2)=w^2$ with
$k\geq 0$. The latter case is excluded, because all $\phi$-fibers are
reduced. Furthermore, if $k$ is even, the closure of $D':=\{x=0\}\cap
\phi^{-1}(U_0)$ is a divisor intersecting the generic fiber in a fiber
of the ruling: a contradiction to the homology
lemma~\ref{homology_lemma} or to $D'$ being $\QZ$-\name{Cartier}. The
remaining case occurs indeed, as was shown in section~\ref{V3ex}.
\end{proof}
\begin{proof}[Step 4: End of the Proof]
After a similar argumentation for the part of $X$ over $U_\infty
=\Pone\setminus\{0\}$, we again obtain the equations of one of the
quadric-degenerations described in section~\ref{V3ex}. Note that the
transition map must commute with the action of $SL_2$. On the other
hand, the only automorphisms of the smooth quadric commuting with the
diagonal action of $SL_2$ are the identity and the involution which
interchanges the horizontal and vertical directions. But
$H^1(\Pone,\ZZ_2)=0$, so that either choice gives a variety which is
isomorphic to one of the examples.
\end{proof}

\subsection{The Construction of the $\mathbb F_4$-Degenerations}\label{V4ex}
Now we consider the case where $X_\eta\cong \Ptwo$. 
In analogy with the construction of the quadric degenerations, set
$V_0 := \Pfive\times \CZ$ with coordinates $([a:b:c:e:f:g],\lambda)$
and let $SL_2$ act on $V_0$ via it's 5-dimensional irreducible
representation on $a\ldots f$. For a given $k\in \NZ$, let the group
$H^*\cong \CZ^*$ act on $V_0$ by 
$$
\xi:([a:b:c:e:f:g],\lambda)\mapsto
([a:b:c:e:f:\xi^{-2k}g],\xi^2\lambda).  
$$
and define $X^k_{0,q}$ to be the variety given by the ideal
\begin{align*}
3e^2-8cf+4f\lambda^kg, 
&&ce-6bf+e\lambda^kg, \\
3be-48af+2c\lambda^kg+2(\lambda^kg)^2,
&&c^2-36af+2c\lambda^k g+(\lambda^k g)^2, \\
bc-6ae+b\lambda^k g, 
&&3b^2-8ac+4a\lambda^k g.
\end{align*}
Note that for a given $\lambda\in \CZ^*\subset \Pone$, the fiber
$X_\lambda$ is isomorphic to $\Ptwo$; the embedding is given by
$[x:y:z]\rightarrow 
([x^2:2xy:2xz+y^2:2yz:z^2:\lambda^{-k}(4xz-y^2)],\lambda)$.

Given another number $l\in \NZ$, construct $X^l_{\infty,q}\subset
V_\infty = \Pfive\times\CZ$ with $H^*$-action given by
$\xi:([a:b:c:e:f:g],\lambda)\rightarrow
([a:b:c:e:f:\xi^{2l}g],\xi^{-2}\lambda)$. The same calculations as in
section~\ref{V3ex} show that $X^k_{0,q}$ and $X^l_{\infty,q}$ glue
together to a variety $X^{(k,l,q)}$ via the relation 
\begin{eqnarray*}
& ([a_0:b_0:c_0:e_0:f_0:g_0],\lambda_0)\sim
([a_\infty:b_\infty:c_\infty:e_\infty:f_\infty:g_\infty],\lambda_\infty)\\
& :\Leftrightarrow \lambda_0\lambda_\infty=1 \text{\ and\ }
[a_0:b_0:c_0:e_0:f_0:g_0] =
[a_\infty:b_\infty:c_\infty:e_\infty:f_\infty:g_\infty\lambda_\infty^{k+l}]. 
\end{eqnarray*}

It is still to be shown that $X^{(k,l,q)}$ has $\QZ$-factorial terminal
singularities and it suffices to show this
for $X^k_{0,q}$. Define $X^k_0$ as in section~\ref{V3ex}, even if $k$
is not odd. Let $\ZZ_2$ act on $X^k_0$ by 
$$
(-1):([x:y:z:w],\lambda)\rightarrow ([x:y:z:-w],\lambda) 
$$ 
We claim that $X^k_{0,q}$ is the quotient of $X^k_0$ by $\ZZ_2$. The
quotient map is given by 
$$
 ([x:y:z:w],\lambda)\rightarrow
([x^2,2xy:2xz+y^2:2yz:z^2:w^2],\lambda).  
$$
and a direct calculation shows that the quotient is isomorphic to
$X^k_{0,q}$.  See \cite[p.~391]{Reid87} for the fact that the
singularities of the quotient are terminal. In order to show that they
are $\QZ$-factorial, it is sufficient to see that all
$\ZZ_2$-invariant divisors in $X^k_0$ are $\QZ$-factorial, if
restricted to the quasi-projective parts. If $k$ is odd, this was
shown for any divisor. If $k$ is even and $D\subset X^k_{0,q}$ a
$\ZZ_2$-invariant divisor, one can argue similarly and use the fact
that
$$
D':=\sum_{\xi\in H^*, \xi^{2k}=1}\xi D 
$$
is a multiple of $D$ and \name{Cartier}.

The same argumentation as in section~\ref{V3ex} shows that
$X^{(k,l,q)}$ can be \name{Mori}-contracted to $\Pone$.

\subsection{The Characterization of the $\mathbb F_4$-Degenerations} 
This is in full analogy to the quadric case.
\begin{prop}[Characterization of the $\mathbb F_4$-Degenerations]
Under the assumptions of~\ref{curve:assumption}, if the generic
$\phi$-fiber is isomorphic to $\Ptwo$, then $X$ is either a linear
$\Ptwo$-bundle or one of the $\mathbb F_4$-degenerations constructed
in section~\ref{V4ex}.
\end{prop}
\begin{proof}
If $X$ is smooth, take a one-parameter subgroup $H<G$ acting
non-trivially on the base $Y$. Given a generic fiber $X_\eta$, there
will always be a line $L\subset X_\eta$, invariant under the action of
the isotropy group $H_\eta$. Then $D:=\overline{H.L}$ is a relatively
ample divisor intersecting $X_\eta$ in $L$. See
\cite[lem.~2.12]{Fuji85} for the fact that this yields an embedding of
$X$ into $\mathbb{P}(\phi_*\mathcal{O}(D))$ which is a
$\Ptwo$-bundle. This must be an isomorphism. Note that $X$ is
automatically smooth if $S$, the semi-simple part of $G$, acts
non-trivially on $Y$.

If $X$ is singular and there is a subgroup $S' < S$, $S'\cong SL_2$,
acting trivially on $Y$ and having a fixed point on generic fibers,
then the subvariety $\{x\in X|\dim S'.x <2\}$ contains a divisor $D$
which intersects $X_\eta$ in an $S'$-homogeneous line. Now argue as in
the proof of proposition~\ref{PropV3ex}. Note that, since $D$ does not
contain a fixed point, it is \name{Cartier}.

It remains to consider the case where $S\cong SL_2$ acts trivially on
$Y$ and stabilizes a quadric curve in $X_\eta$. As above, let $D$ be
the union of these curves. In complete analogy to the proof of
proposition~\ref{PropV3ex}, all fibers are isomorphic to $\Ptwo$ or
$\mathbb F_4$, $D$ is \name{Cartier} and yields an embedding into a
$\Pfive$-Bundle $\mathbb{P}(E)$. Here $E$ splits $S$-equivariantly
into a direct sum of a 5-dimensional bundle $E_5$, where $S$ acts via
it's irreducible representation, and a 1-dimensional bundle $E_1$
where the $S$-action is trivial. Furthermore, the subbundle
$\mathbb{P}(E_5)$ is the unique hyperplane intersecting $X$ in $D$!

We continue to argue as in~\ref{PropV3ex}, using the fact that all
$SL_2$-invariant subsets in $\Pfive$, isomorphic to $\Ptwo$ and not
contained in the $SL_2$-invariant hyperplane are given by
\begin{align*}
3e^2-8cf+4f\lambda g, 
&& ce-6bf+e\lambda g, \\
3be-48af+2c\lambda g+2\lambda^2g^2, 
&& c^2-36af+2c\lambda g+(\lambda g)^2, \\
bc-6ac+b\lambda g,
&& 3b^2-8ac+4a\lambda g,
\end{align*} 
where $\lambda\in \CZ^*$. Consequently, $X$ is locally given by the
equations from section~\ref{V4ex}. There is no choice of how the
affine parts can be $SL_2$-equivariantly glued.
\end{proof}

\section{Minimal Models over surfaces}\label{surfsect}

The primary aim of this section is to classify the relatively minimal
varieties over surfaces. The following lemma describes a particularly
simple situation.

\begin{lem}\label{Surf:S_acts_triv_on_Y}
In the situation of theorem~\ref{mainthm}, let $Y$ be a surface. If
$S<G$ is a semi-simple group which acts trivially on $Y$, then $X\cong
Y\times \Pone$. In particular, $X$ and $Y$ are smooth
\end{lem}
\begin{proof}
Since $\phi$ is an extremal contraction, all fibers must be of
dimension 1. 

Note that $S$ acts transitively on the generic fibers. Thus, $S\cong
SL_2$ and $S$ has no fixed points: a linearization of the
$SL_2$-action would give a contradiction.

Consequently, $S$ acts transitively on all fibers, and if $U<S$ is a
maximal connected unipotent subgroup, and $\Sigma$ it's fixed point
set, then $X=S.\Sigma\cong \Pone\times \Sigma$. 
\end{proof}

Due to the preceding lemma, we may consider for the rest of this
section that the semi-simple part of $G$ acts non-trivially on $Y$:
\begin{assumption}\label{surf:assumption}
Let $X$ and $G$ be as in theorem~\ref{mainthm} and let
$\phi:X\rightarrow Y$ be an extremal contraction to a surface. Let
$S<G$ be a maximal semi-simple subgroup and assume that no simple
factor of $S$ acts trivially on $Y$.
\end{assumption}

\subsection{Models over Singular Surfaces}
We start with the construction of the relatively minimal varieties
over a singular surface.
\begin{example}\label{Surf:SingSurf}
Set $\tilde{X} := \mathbb{P}(\mathcal{O}_{\Ptwo}(e) \oplus
\mathcal{O}_{\Ptwo})$. The automorphism group of $\tilde{X}$ contains
a product $G:=SL_2 \times \CZ^*$, where $SL_2$ has a fixed point in
$\Ptwo$, 2-dimensional orbits in $\tilde{X}$ and acts trivially on the
fiber over the fixed point. The factor $\CZ^*$ acts in fiber direction
only, i.e.~trivially on $\Ptwo$.  Embed $\ZZ_2$ diagonally into $G$,
i.e.~consider the subgroup generated by $(Diag(-1,-1),-1)$. Then
$X:=\tilde{X}/\ZZ_2$ is a singular model over $\mathbb{F}_2$, the cone
over the rational normal curve in $\Ptwo$.
\end{example}

We will see that these are the only possibilities.

\begin{notation}
Call a divisor $D\subset X$ a ``rational section'' iff it intersects
the generic $\phi$-fiber with multiplicity 1. Note that a rational
section is a section iff it does not contain a whole $\phi$-fiber.
\end{notation}

\begin{prop}
Under the assumptions~\ref{surf:assumption}, let $Y$ be singular. Then
$X$ is one of the varieties constructed in example~\ref{Surf:SingSurf}.
\end{prop}
\begin{proof}
As a first step, construct a rational section. By assumption, $S$ acts
non-trivially on $Y$. It follows from the classification that $Y \cong
\mathbb{F}_n$, the cone over a rational normal curve and $S \cong
SL_2$. The $S$-isotropy $S_\eta$ of a generic point $\eta\in Y$ is an
extension of a maximal unipotent group by a cyclic group. Thus,
$S_\eta$ fixes at least one point in the fiber $X_\eta$ so that the
closure $E$ of at least one $S$-orbit is a rational section indeed.

The \name{Weil}-divisor $E$ is not \name{Cartier}, or else use
\cite[lem.~2.12]{Fuji85} and obtain a contradiction to ``$Y$
singular''. Thus, $X$ is singular.

The next step is to construct a cover of $X$. Observe that a fiber
$X_\mu$ through the singular set $\Sing(X)$ is pointwise $S$-fixed and
linearize the $S$-action at a generic point $f
\in X_\mu$. After proper choice of coordinates, one may identify a
neighborhood $U(f)\cong \Delta_1\times \Delta_2$, where $\Delta_1$ is
a one-dimensional and $\Delta_2$ a 2-dimensional ball. We can assume
that $S$ acts only on the second component and that the map
$\phi|_{U(f)}$ is given by the projection to the second factor
followed by taking the quotient by $\ZZ_n$.

Let $\gamma:\Ptwo\rightarrow \mathbb{F}_n$ be the natural cyclic
cover.  Observe that $\gamma$ is $S$-equivariant and set $X':= X
\times_{\mathbb{F}_n} \Ptwo$. Calculating the preimage of $U(f)$ one
obtains $\Delta_1\times(\text{$n$ copies of\ }\Delta_2 \subset \CZ^4
\text{\ meeting in a point})$. If $\tilde{X}$ is the normalization of
$X'$, the preimage of $U(f)$ becomes $\Delta_1 \times (\Delta_2\coprod
\ldots \coprod \Delta_2)$, so that $\tilde{X}$ is a $n$:1 cover over
$X$, with finite singular set. The calculation also shows that
$\tilde{X}$ is \name{Galois} with group $\Gamma = \ZZ_n$.

We claim that $\tilde{X}$ is a split linear $\Pone$-bundle over
$\Ptwo$: $\tilde{X} \cong \mathbb{P}(\mathcal{L}\oplus\mathcal{O})$. 
If $\tilde{\phi}:\tilde{X}\rightarrow \Ptwo$ is the natural map,
consider the $\tilde{\phi}$-fiber $\tilde{X_\mu}$ over the unique 
$S$-fixed point in $\Ptwo$. As it's image in $X$ is pointwise
$S$-fixed, $\tilde{X_\mu}$ is, too. Using the linearization
argument, let $\tilde{E}$ be the closure of a generic $S$-orbit,
intersecting $\tilde{X_\mu}$ in a generic point. Observe that
$\tilde{E}$ contains a unique $S$-fixed point and is smooth
there. Consequently, $\tilde{E}$ is a smooth section, $\tilde{X} =
\mathbb{P}(\mathcal{E})$ is a linear $\Pone$-bundle and, as all
Ext-groups on $\Ptwo$ vanish, $\mathcal{E}$ is split: we may assume
$\mathcal{E} = \mathcal{O}_{\Ptwo}(e) \oplus \mathcal{O}_{\Ptwo}$. 

We must show that the action of $\Gamma$ is the same as in the
example~\ref{Surf:SingSurf} above. Identify an $SL_2$-invariant
neighborhood of $\tilde{X_\mu}$ with $\CZ^2 \times \Pone$ in a way
that $S$ acts on the first factor only. By equivariance, $\Gamma$ maps
$S$-orbits to $S$-orbits. Consequently, the quotient by $\Gamma$ has
two cyclic quotient singularities of type $\frac{1}{n}(1,1,a)$ and
$\frac{1}{n}(1,1,-a)$. As quotient singularities are terminal only if
of type $\frac{1}{n}(1,a,-a)$ (cf. \cite[sect. 5.3]{Reid87}), $n=2$
and $a=1$.  This yields the claim.
\end{proof}

\subsection{Models over Smooth Surfaces}
In \cite{K98a} we were discussing the possibility to compactify
homogeneous spaces to particularly simple varieties. We refer the
reader to section~5.2 of that paper for a proof of the fact that the
$X$ is a automatically a linear $\Pone$-bundle if $X$ is relatively
minimal over a smooth surface and $G$ is not solvable. The rest of
section~\ref{surfsect} is concerned with an investigation which
rank-two vector bundles do actually occur. We assume that $X$ is a
linear bundle without further mention.

Remark that under the assumptions~\ref{surf:assumption} $Y$ is a
rational $G$-quasihomogeneous surface with non-trivial
$S$-action. Since $X$ is now supposed to be smooth, $Y$ is smooth, so
that $Y\cong \Sigma_n$ or $\Ptwo$. Later on, we will consider these
cases separately.

\begin{notation}\label{Fnot}
Let $\phi:X\rightarrow Y$ be as above and assume that
there exists a map $\pi: Y\rightarrow Z\cong\Pone$, e.g. if $Y$ is
isomorphic to a (blown-up) \name{Hirzebruch} surface $\Sigma_n$. Then,
if $F \in Z$ is a generic point, set $F_Y:=\pi^{-1}(F)$ and
$F_X:=\phi^{-1}(F_Y)$. 

\end{notation}

\subsection{The Construction of the Diagonally Twisted Bundles}
The following varieties will be of great importance in the classification:

\subsubsection{The Construction of the $X_{\Sigma_n,k_0,k_\infty}$}
\label{ex01}
Let $Y$ be the \name{Hirzebruch}-surface $\Sigma_n$, $n >0$ and $X:= Y
\times \Pone$. Let $S:=SL_2$ act on $Y$ and $\Pone$ and let $S$ act on
$X$ diagonally, i.e.~simultaneously on both components.

We claim that $S$ acts quasihomogeneously on $X$ and that the
exceptional set (i.e.~the complement of the open orbit) contains a
unique $S$-invariant section over $Y$. In order to see this, let $B<S$
be the \name{Borel} group of $S$ stabilizing $F_X$. The $B$-action on
$F_X$ is very special: Since the $S$-action on $\Sigma_n$ stabilizes
the 0- and the $\infty$-sections, the $B$-action on $F_X$ stabilizes
two fibers. Therefore the unipotent part $B_U$ of $B$ acts in fiber
direction only, showing that the $B$-action on $F_X$ is
quasihomogeneous and that there is exactly one $B$-invariant section
in $F_X$. Using the $S$-action in order to move $F_X$ around shows
that $S$ does indeed act quasihomogeneously and that there is a unique
$S$-invariant rational section $E$. The fact that $S$ does not have
any fixed points on $Y$ immediately implies that $E$ is indeed a
section and $E\cong Y \cong \Sigma_n$. Let $E_0$ and $E_\infty$ denote
the 0- and $\infty$-section of $E$, respectively.

The curves $E_0$ and $E_\infty$ are the only $S$-invariant subsets in
$E$. We can now perform an elementary transformation with center being
$E_0$ or $E_\infty$, obtaining a new $\Pone$-bundle which is not
necessarily the compactification of a line bundle. By elementary
transformation we understand the process of blowing up $E_0$ and then
blowing down the strict transform of $\phi^{-1}\phi(E_0)$. Such
transformations always exist; see \cite{Maru73} for a complete
reference. Since the centers of the transformations are
$SL_2$-invariant, $SL_2$ acts on the transformed varieties, and the
entire procedure is equivariant.

The strict transform of $E$ is again invariant and isomorphic to
$\Sigma_n$ so that one may iterate the process. Let
$X_{\Sigma_n,k_0,k_\infty}$ be the variety obtained by transforming
$k_0$ times with center being the 0- and $k_\infty$ times with center
being the $\infty$-section of $E$. Let $F_{k_0,k_\infty}$ be the
strict transform of $F_X$ in $X_{\Sigma_n,k_0,k_\infty}$.

As above, $B_U$ acts on $F_{k_0,k_\infty}$ by adding multiples of a
section. Note that $F_{k_0,k_\infty} \cong \Sigma_{k_0+k_\infty}$ and
the sections added by $B_U$ vanish of order $k_0$ at $F_{k_0,k_\infty}
\cap E_0$ and of order $k_\infty$ at $F_{k_0,k_\infty} \cap E_\infty$.

\subsubsection{The Construction of the $X_{\Sigma_0,n}$}
\label{Surf:X_Sigma_0_n}
Let $S:=SL_2$ act diagonally on $Y=\Sigma _0$. Since $S$ is a
simply-connected semi-simple group and $H^1(Y,\mathcal O)=0$, the
$S$-action on $Y$ can be lifted to the total space of any line bundle
${\mathcal O}(n,m)$ over $Y$; see \cite[p.~98]{H-Oe} for details. For
$n\in \NZ^+$ the group $S$ therefore acts on the compactification
$X=\mathbb P({\mathcal O}(n,-n)\oplus {\mathcal O})$ which is a
$\Pone$-bundle $\phi :X\to Y$.  This lifting is unique up to the
$\CZ^*$-action given by the principal $\CZ^*$-actions on the first
factor.

Let $\sigma _0$ and $\sigma _\infty$ be the $S$-invariant sections
defined by the direct sum structure. Since there are no other
sections, it follows that $S$ acts transitively on the complement
$X\setminus (\sigma _0\cup \sigma _\infty \cup \Delta _X$), where
$\Delta _X$ is the preimage $\pi ^{-1}(\Delta )$ of the $S$-invariant
diagonal $\Delta $ in $Y$.

If $i:\Delta \hookrightarrow Y$ is the canonical embedding, then
$i^*({\mathcal O}(n,-n))$ is trivial.  Thus $\phi |_{\Delta _X}:\Delta
_X\to Y$ is the trivial $\Pone$-bundle and therefore all $S$-orbits
in $\Delta _X$ are $1$-dimensional sections over $\Delta =\Pone$.

Let $C=Sx$ be such a section which does not lie in $\sigma _0\cup
\sigma _\infty $ and define $X_{\Sigma _0,n}$ to be the elementary
transformation of $X$ with respect to $C$ in $\Delta _X$.  This
manifold is still an $S$-equivariant $\Pone$-bundle over
$Y$. However now the transforms $\sigma _0'$ and $\sigma _\infty '$
intersect transversally in an $S$-orbit $C'=Sx\cong \Pone$ over
$\Delta $.

Given any two $S$-orbits $C_i:=Sx_i$ as above, there exists a unique
transformation $g$ of the $\CZ^*$-action which commutes with the
$S$-action so that $g(C_1)=C_2$. This defines an $S$-equivariant
isomorphism between the spaces $X^1_{\Sigma _0,n}$ and $X^2_{\Sigma
_0,n}$ which are defined by elementary transformations along $C_1$ and
$C_2$ respectively.  In this sense the {\em diagonally twisted bundle}
$X_{\Sigma _0,n}$ is uniquely defined.

\subsection{The Classification of $S$-quasihomogenous Bundles}

\subsubsection{Bundles over $\Sigma_n$} The following lemma gives a
first characterization of split linear bundles:
\begin{lem}\label{lemww}
Under the assumptions~\ref{surf:assumption}, assume additionally that
$Y\cong \Sigma_n$ and that $F_X \cong \Pone\times\Pone$, where $F_X$ is
defined as in notation~\ref{Fnot}. Then $X$ is isomorphic to a fibered
product: $X\cong Y\times_Z Y'$. In particular, $X$ is a split linear
bundle: $X\cong
\mathbb{P}(\mathcal L\oplus \mathcal O)$.
\end{lem}
\begin{proof}
The space $X$ has relative \name{Picard}-number 2 over $Z$ and the
general fiber $F_X$ is \name{Fano}. Thus, there exists a second
\name{Mori}-contraction $\phi':X\rightarrow Y'$ over $Z$ which is
different from $\phi$.

The \name{Picard}-number of $Y'$ is 2, so $Y'$ is not a curve. If
$\dim Y'=3$, then the contraction was divisorial inducing a
contraction from $F_X$ to a surface, which is obviously
impossible. Therefore, the contraction $\phi'$ is of fiber type.
Consequently $X$ is a $\Pone$-bundle over $Y'$ with fibers being the
horizontal curves in $F_X$ and their translates. 
\end{proof}

With the aid of the preceding lemma we can now carry out the
classification of $S$-quasihomogeneous bundles over $\Sigma_n$.

\begin{prop}[Characterization of $X_{\Sigma_n,k_0,k_\infty}$]\label{hhh}
Under the assumptions~\ref{surf:assumption}, if $Y\cong \Sigma_n$,
$n>0$ and $S$ acts almost transitively on $X$, then either $X\cong
Y\times\Pone$, or $S\cong SL_2$ and there exist numbers $k_0, k_\infty
\geq 0$ such that $X\cong X_{\Sigma_n,k_0,k_\infty}$ is one of the
diagonally twisted bundles constructed in section~\ref{ex01}.
\end{prop}
\begin{proof}
Since no factor of $S$ acts trivially on $Y$, $S\cong SL_2$. There are
exactly two $S$-invariant curves $\sigma_0$, $\sigma_\infty$ in $Y$;
these are sections over $Z$. Call the preimages $\phi^{-1}(\sigma_0)$
and $\phi^{-1}(\sigma_\infty)$ of these sections $A_0$ and $A_\infty$,
respectively.

{\em The $S$-invariant divisors in $X$:} Let $F_X$ be as in
notation~\ref{Fnot}. If $B$ is the \name{Borel} group in $SL_2$
stabilizing $F_X$, then, because $F_X \cap A_0$ and $F_X \cap
A_\infty$ are invariant, $B_U$, the unipotent part of $B$, acts
trivially on the base. Instead, $B_U$ acts on the fibers of
$\phi|_{F_X}$ and fixes a unique point in each. Consequently there
exists a unique $B$-invariant section in $F_X$; other $B$-invariant
curves are the fibers $A_0\cap F_X$ and $A_\infty\cap F_X$. Using $S$
to move $F_X$, one sees that the only closed $S$-invariant divisors
are $A_0$, $A_\infty$ and a unique section, called $E$. Furthermore,
$E\cap F_X$ being the only $B$-invariant section implies that $E\cap
F_X$ is the unique curve of negative self-intersection in $F_X$ if
$F_X
\cong \Sigma_m$, $m>0$.

{\em Application of the Algorithm:} As a next step perform the
sequence of elementary transformations given by the algorithm outlined
in figure~\ref{algtab}. One must show that the algorithm stops. Since
the center of the elementary transform intersects $F_X$ in a point not
contained in $E$, i.e.~not contained in the $\infty$-section of $F_X$,
the self-intersection of $E\cap F_X$ in $F_X$ rises by one. Since it
was negative when the algorithm started, it will eventually become
zero, implying $F_X\cong \Sigma_0$, and the process terminates.
\begin{figure}
$$
\xymatrix{
 & {\txt{Start}} \ar[d] \\
 & {\txt{Is $F \cong \Sigma_0$?}} \ar[d]_{\txt{no}} \ar[r]_{\txt{yes}}
 & {\txt{Stop (A)}}	\\
 {\txt{Perform an elementary\\ transformation
  with center\\ being the curve not \\contained in $E$.}} \ar `u[u][ru]
 & {\txt{Are there two $S$- \\ invariant curves in $A_0$?}} 
   \ar[d]_{\txt{no}} \ar[l]^{\txt{yes}}	\\
 & {\txt{Are there two $S$- \\ invariant curves in $A_\infty$?}}
   \ar[d]_{\txt{no}} \ar `l[l]^{\txt{yes}} [lu]	\\
 & {\txt{Stop (B)}} \\
}
$$
\caption{an algorithm for simplifying special $\Pone$-bundles
\label{algtab}}
\end{figure}

We claim that the point ``Stop (B)'' is never arrived at, i.e.~$A_0$
and $A_\infty$ having only one $S$-invariant curve implies $F_X\cong
\Sigma_0$. Note that $S$ has only one invariant curve in $A_0$ and 
$A_\infty$ if and only if $A_0$, $A_\infty \cong \Sigma_0$ and $S$
acts diagonally. This implies that $B_U$ has unique fixed points in
$F_X \cap A_0$ and $F_X \cap A_\infty$, namely the intersection with
$E$. Consequence: if $\sigma\subset F_X$ is a section not intersecting
$E$ and $u\in B_U\setminus \{1\}$, then $\sigma$, $u.\sigma$ and
$E\cap F_X$ are three mutually disjoint sections in $F_X$ over $F_Y$!
Therefore $F_X \cong \Sigma_0$.

{\em The Situation where the Algorithm stops:} Let us now assume that
the algorithm already stopped, i.e.~$F_X\cong \Sigma_0$. Apply
lemma~\ref{lemww}: as the algorithm terminates, the transformed
variety is isomorphic to $Y\times_Z Y'$. Now to say that there is a
unique $B$-invariant section in $F_X$ over $F_Y$ which is not diagonal,
it is equivalent to say that there exists a unique $S$-invariant curve in
$Y'$. Hence $Y' \cong \Sigma_0$ and $SL_2$ acts diagonally. In
particular, $X$ is the trivial bundle over $Y$ and $SL_2$ acts
diagonally. Recall that this is the starting situation of
section~\ref{ex01}.

{\em End of the Proof:} As a last step there is to prove that the
inverses of the transformations we performed are the transformations
used in section~\ref{ex01}, i.e.~elementary transformations with
center being $E_0$ or $E_\infty$. This, however, is clear if one takes
into account that the algorithm transforms with centers being curves
in $A_0$ or $A_\infty$ {\em not} intersecting $E$.
\end{proof}

\subsubsection{Bundles over $\Sigma_0$} The primary goal of this
section is to characterize the diagonally twisted $\Pone$ bundles over
$\Sigma_0\cong \Pone\times\Pone$. It is necessary to prove that
sections which arise as closures of $S$-orbits are either disjoint or
intersect transversally. The following lemma is a first step in this
direction. 
\begin{lem}\label{Surf:Curves_in_Sigma_n}
Let $B<SL_2$ be a \name{Borel} group and $\Sigma_n$ a
\name{Hirzebruch}-surface with a surjection $\phi:\Sigma_n \rightarrow
\Pone$. Assume that $B$ acts almost transitively on $\Sigma_n$ and
that $\Sigma_n$ contains two $B$-invariant sections $\sigma_1$ and
$\sigma_2$ over $\Pone$. Then either $\sigma_1$ and $\sigma_2$ are
disjoint or they intersect transversally.
\end{lem}
\begin{proof}
As a first step, remark that there are at most 2 $B$-invariant sections
in $\Sigma_n$. The existence of a third would contradict the almost
transitive action, because if $\eta\in \Pone$ is a point in the open
orbit, then it's isotropy group must act almost transitively on the
fiber $X_\eta$ and fix the intersection with of $X_\eta$ all invariant
sections. But there are no non-trivial automorphisms of a generic
fiber fixing three points. This means that we only need to find two
disjoint or transversal sections in order to prove the claim.

The same line of argument shows that $B$ may not have two fixed points
on the base $\Pone$, for otherwise the unipotent part $U$ of $B$ would
act trivially on $\Pone$. Thus, if $X_\eta$ is a general fiber, $U$
would stabilize $X_\eta$. But $U$ acts non-trivially on $\Sigma_n$,
because $\Sigma_n$ and $B$ both have dimension 2. Consequently, $U$
would act non-trivially on $X_\eta$, and $X_\eta$ contains only one
point which is invariant under the isotropy group $B_\eta$. So there
would only be one invariant section.

We prove the lemma by induction

{\em Start of Induction: $n=0$.} Assume without loss of generality
that $\phi:\Sigma_0\cong \Pone\times\Pone\rightarrow \Pone$ is the
projection onto the first factor. If the $B$-action on the second
factor has two fixed points, then there are 2 disjoint sections, and
we are finished. Otherwise, note that there is only one $B$-action on
$\Pone$ with exactly one fixed point ---up to isomorphy. Thus, after
appropriate choice of coordinates, we can assume that $B$ acts
diagonally on $\Sigma_0$. In this situation $B$ stabilizes the
diagonal in $\Sigma_0$ and a fiber of the projection to the second
factor. These curves meet transversally.

{\em Step of Induction:} Assume that the lemma is true for all numbers
smaller than a given $n>0$. We will assume that the lemma is false for
$n$ and derive a contradiction. Thus, suppose that we are given two
$B$-invariant divisors $\sigma_1$ and $\sigma_2$ which do not intersect
transversally. Let $\sigma_1$ be the unique curve of negative
self-intersection in $\Sigma_n$, this curve is a section which is
invariant under the full automorphism group. Let $F$ be the unique
$B$-invariant $\phi$-fiber, the preimage of the $B$-fixed point in
$\Pone$.

{\em Claim:} the group $B$ has two fixed points in $F$. 

If the claim holds, then we can perform a $B$-equivariant elementary
transformation where we choose the center to be the $B$-fixed point
which is {\em not} contained in $\sigma_1$. If $X'$ is the transformed
variety and $\sigma'_1$ and $\sigma'_2$ are the strict transforms of
$\sigma_1$ and $\sigma_2$, then $\sigma'_1$ and $\sigma'_2$ still
intersect non-transversally: the intersection number
$\sigma'_1.\sigma'_2$ is even bigger than $\sigma_1.\sigma_2$. On the
other hand, by choice of the center, $X'\cong \Sigma_{n-1}$. We obtain
a contradiction to the induction hypothesis and are finished.

It remains to show the claim. Again assume to the contrary,
i.e.~assume that there was only one $B$-fixed point in $F$. Let $T<B$
be a torus. Since all $B$-actions on $\Pone$ which have only one fixed
point are isomorphic, we know that $T$ acts on $F\cong \Pone$ with
weight 2. Similarly, $T$ acts on $\sigma_1$ with weight 2; this is
because $\sigma_1$ is a section and the restricted map
$\phi|_{\sigma_1}:\sigma_1\rightarrow \Pone$ is equivariant. Now
linearize the $T$-action at the intersection point $\sigma_1\cap
F$. Realize that $F$ and $\sigma_1$ intersect transversally. But the
only 2-dimensional $T$-representation space containing two $T$-invariant
curves of weight 2 which additionally intersect transversally has
weights $(2,2)$. Thus, any two $T$-invariant curves passing through the
intersection point must intersect transversally. But since the
intersection $\sigma_1\cap\sigma_2$ is $B$-fixed,
$\sigma_1\cap\sigma_2\subset F$ so that $\sigma_1$ and $\sigma_2$ are
two $T$-invariant curves passing through $\sigma_1\cap F$! This is
absurd.
\end{proof}

\begin{prop}[Characterization of the $X_{\Sigma_0,n}$]
\label{Surf:Char_X_Sigma0_n}
Under the assumptions~\ref{surf:assumption}, if $Y\cong \Sigma_0$ and
$S$ acts almost transitively on $X$, then $X$ is a splitting bundle or
one of the diagonally twisted bundles $X_{\Sigma_0,n}$ from
example~\ref{Surf:X_Sigma_0_n}.
\end{prop}
\begin{proof}
To start with, choose a morphism $\pi:Y\rightarrow Z\times\Pone$ and
define $F$, $F_X$ and $F_Y$ as in notation~\ref{Fnot}. If $F_X\cong
\Sigma_0$, we are finished by using lemma~\ref{lemww}. Thus, assume that
$F_X\cong \Sigma_n$, $n>0$. No simple factor of $S$ acts trivially
on $Y$. Thus either $S\cong SL_2$, acting diagonally on $Y$ or $S\cong
SL_2\times SL_2$.

If $S\cong SL_2\times SL_2$, let $S'<S$ be the factor of $S$ acting
trivially on $Z$ and note that there are always two disjoint
$S'$-invariant sections $\sigma_1$ and $\sigma_2$ in $F_X$ over $F_Y$. If
$S''$ is the other factor of $S$, then $S''$ must stabilize the locus
where $S'$ has 1-dimensional orbits; this is because $S'$ and $S''$
commute. In particular, $S''.\sigma_1$ and $S''.\sigma_2\subset F_X$
are two disjoint sections over $Z$, displaying $X$ as a splitting
linear bundle.

For the remainder of the proof consider the situation where $S\cong
SL_2$. The isotropy $S_\eta$ at a generic point $\eta\in Y$ is a
torus. This torus fixes two points in the associated $S_\eta$-invariant
$\phi$-fiber $X_\eta$, and a standard argument shows that the closures
of their $S$-orbits are rational sections. If these are disjoint, we
can stop here. Thus, assume that they have non-trivial
intersection. We claim that $X\cong X_{\Sigma_0,n}$.

As a first step in this direction show that the $\sigma_\bullet$
intersect transversally. In order to see this, consider the stabilizer
$B$ of $F_X$, which is a \name{Borel}-subgroup of $S$.  The curves
$\sigma_1\cap F_X$ and $\sigma_2\cap F_X$ are two $B$-invariant
sections in $F_X$ over $F_Y$ which intersect in a single
point. Lemma~\ref{Surf:Curves_in_Sigma_n} claims that the intersection
of these curves must be transversal. This transversality implies that
the sections become disjoint if one performs an elementary
transformation with center $\sigma_1\cap \sigma_2$. In other words, if
$X'$ is the transformed variety, then the strict transforms of
$\sigma_1$ and $\sigma_2$ are disjoint. This already shows that $X'$
is a splitting linear bundle.

{\em The Triviality of the Bundle over the Diagonal:} If $\Delta_X$
denotes the preimage of the $S$-invariant diagonal $\Delta\subset Y$,
then $\Delta_X$ contains the center of the transformation and is
transversal to both $\sigma_1$ and $\sigma_2$. Thus, after blowing up
$\sigma_1\cap \sigma_2$, the strict transform of $\Delta_X$ becomes
disjoint from the strict transforms of the $\sigma_\bullet$. This in
turn implies that the exceptional divisor of the blow-up is isomorphic
to $\Sigma_0$, and $S$ acts with one-dimensional orbits there. By
construction, the same holds for the preimage of $\Delta$ in $X'$. So
$X'$ is of the form $\mathbb P(\mathcal O(n,-n)\oplus
\mathcal O)$ since these are the only split $\Pone$-bundles which are
trivial over the diagonal. 

We have seen that the center of the back-transformation $X'\dasharrow
X$ is not contained in one of the $S$-invariant sections. So that
back-transformation is exactly the construction performed in
example~\ref{Surf:X_Sigma_0_n}. 

This establishes the isomorphy $X\cong X_{\Sigma_0,n}$ once we know
that $n\not =0$. Recall that $SL_2$ acts almost transitively on
$X'$. But if $n$ was 0, then $X'\cong \Sigma_0\times\Pone$ was the
trivial bundle and $SL_2$ having one-dimensional orbits over the
diagonal would imply that $SL_2$ acts trivially on the second factor,
a contradiction.

This proves $X\cong X_{\Sigma_0,n}$, and the claim is shown.
\end{proof}

\subsubsection{Bundles over $\Ptwo$} The classification of bundles
over $\Ptwo$ is mainly a corollary of the classifications we have
carried out already.
\begin{prop}\label{Surf:QuasiHom_over_Ptwo}
Under the assumptions~\ref{surf:assumption}, if $Y\cong \Ptwo$ and 
$S$ acts almost transitively on $X$, then $X$ is isomorphic to a
splitting $\Pone$-bundle, to the flag manifold $F_{(1,2)}(3)$, a
blow-down of $X_{\Sigma_1,k_0,0}$, or to a quotient of one of the
models over $\Sigma_0$.
\end{prop}
\begin{proof}
We tell between the possible $S$-actions on $Y$:
\begin{description}
\item[$\boldsymbol{S\cong SL_3}$] If $SL_3$ acts transitively on $X$,
then $X\cong F_{(1,2)}(3)$; this follows from the classification of
the homogeneous manifolds. See \cite{Winkelmann}. Otherwise, the
exceptional set $E$ is an unbranched cover of $\Ptwo$. A connected
component of $E$ is a section, realizing $Y$ as a splitting linear
bundle.

\item[$\boldsymbol{S\cong SL_2}$, and $\boldsymbol S$ has a fixed
point $\boldsymbol{\mu \in Y}$] Blow up the $\phi$-fiber $X_\mu$ an
obtain a $\Pone$-bundle $X'\rightarrow Y'\cong \Sigma_1$. Let
$E\subset X'$ be the exceptional divisor of the blow-up. By
proposition~\ref{hhh}, there are only two possibilities:

If $X'$ is a splitting linear bundle, then let $\sigma_1$ and
$\sigma_2$ be two disjoint sections over $Y$. They intersect $E$ in
two different fibers of the fibration $E\rightarrow
X_\mu$. Consequence: the images of $\sigma_1$ and
$\sigma_2$ in $X$ are disjoint sections, too. Thus $X$ is split.

If $X'\cong X_{\Sigma_1,k_0,k_\infty}$, claim that $k_\infty=0$. As a
first step, realize that $S$ acts non-trivially on $X_\mu$, or else a
linearization argument would reveal that $S$ has only 2-dimensional
orbits. Thus, $S$ acts diagonally on $E$, and there is exactly one
$S$-invariant curve in $E$. This already shows that $k_\infty=0$, for
otherwise $E$ had to contain two distinct $S$-invariant curves: the
intersection with the unique $S$-invariant section over $Y'$ and the
center of the back-transformation $X_{\Sigma_1,k_0,k_\infty-1}
\dasharrow X_{\Sigma_1,k_0,k_\infty}$.

\item[$\boldsymbol{S \cong SL_2}$, and $\boldsymbol S$ does not have a
fixed point in $\boldsymbol Y$] Recall that there exists an
$S$-equivariant cover $\gamma: \Sigma_0 \rightarrow \Ptwo$. The
pull-back $X' := X \times_Y \Sigma_0$ is $S$-quasihomogeneous. It
was shown in proposition~\ref{Surf:Char_X_Sigma0_n} that $X'$ is a
splitting bundle or $X'\cong X_{\Sigma_0,n}$, and $X$ is a quotient of
$X'$ by $\ZZ_2$.
\end{description}
\end{proof}

\subsection{The Remaining Cases} It remains to consider the cases
where $S$ does not act almost transitively. We start with a
classification of the models over $\Sigma_n$.

\begin{prop}\label{ggg}
Under the assumptions~\ref{surf:assumption}, if $Y\cong \Sigma_n$ and
$S$ has generic orbits of dimension $\leq 2$, then $X$ is a splitting
linear bundle.
\end{prop}
\begin{proof}
If $S$ acts almost transitively on $Y$, then there is a subgroup
$S'<S$ acting almost transitively on $Y$ with $S'\cong SL_2$. Hence,
assume without loss of generality that $S\cong SL_2$. Choose
$\pi:Y\rightarrow Z\cong \Pone$. If $F_X\cong \Sigma_0$, apply
lemma~\ref{lemww} and stop. Otherwise, let $B<S$ be the
\name{Borel}-group stabilizing $F_X$ and note that the generic
$B$-orbit in $F_X$ has dimension at most 1.

{\em Claim:} there are two disjoint $B$-invariant sections $\sigma_0$,
$\sigma_\infty$ in $F_X$ over $F_Y$. 

In order to prove the claim, let $T<B$ be a maximal torus. Recall that
$T$ is not normal in $B$. Thus $T$ acts non-trivially on $F_X$, and a
curve on $F_X$ is $B$-invariant if and only if it is
$T$-invariant. Now the claim follows from the following fact: a
maximal torus in $\Aut(\Sigma_m)$ always contains a subtorus whose
fixed point set are two sections.

Now we apply the claim: $D_\bullet := S.\sigma_\bullet$ are disjoint
sections in $X$ over $Y$, and we are finished.

If $S$ acts on $Y$ with 1-dimensional orbits, then choose
$\pi:Y\rightarrow Z$ so that $S$ acts trivially on $Y$. Now there are
two possibilities: the first is that $F_X\cong \Sigma_m$, where
$m>0$. But then there are necessarily two disjoint sections in $F_X$
over $F_Y$. Recall that the set $D:=\{x\in X|\dim S.x \leq 1\}$ is
closed. By what we saw above, $D$ is a $2:1$ unbranched cover over
$Y$. But $Y$ is simply connected. Thus, $D$ consists of 2 disjoint
sections, and $X$ is a split linear bundle.
\end{proof}

It remains to consider models over $\Ptwo$.
\begin{prop}
In the setting of~\ref{surf:assumption}, if $Y\cong \Ptwo$, and $S$
has generic orbits of dimension $\leq 2$, then $X$ is a splitting
bundle.
\end{prop}
\begin{proof}
Consider the different possibilities for the $S$-action on $Y$.
\begin{description}
\item[$\boldsymbol{S\cong SL_3}$] All $SL_3$-orbits are isomorphic to
$\Ptwo$ and, by $S$ acting transitively, are unbranched covers of
$Y$. Three of them yield the identification
with the trivial bundle.

\item[$\boldsymbol{S\cong SL_2}$, and $\boldsymbol S$ has a fixed
point $\boldsymbol{\mu \in Y}$] We blow up the $\phi$-fiber $X_\mu$,
obtain a $\Pone$-bundle over $\Sigma_1$ and apply
proposition~\ref{ggg}. Argue as in the proof of
proposition~\ref{Surf:QuasiHom_over_Ptwo} to see that $X$ is split as
well. 

\item[$\boldsymbol{S \cong SL_2}$, and $\boldsymbol S$ does not have a
fixed point in $\boldsymbol Y$] Take a \name{Borel} group $B<S$. The
isotropy $B_\eta$ of a generic point in $\eta \in Y$ is finite and
cyclic. Hence there exists unique $B_\eta$-invariant point $f\in
X_\eta$ and $D:=\overline{B.f}$ is a unique $S$-invariant section. The
vanishing of all Ext-groups yields the claim.
\end{description}
\end{proof}

\section{Minimal Models over a point}

Now we classify the situation of theorem~\ref{mainthm} under the
additional assumption that $Y$ is a point. The next lemma shows that
nontrivial models occur only if the semi-simple part $S$ of $G$ is
isomorphic to $SL_2$.
\begin{lem}
Under the assumptions of theorem~\ref{mainthm}, let $Y$ be a point. If
$G$ contains a semi-simple group other than $SL_2$, then $X\cong
\Pthree$, $\QZ_3$ or the weighted projective space 
$\mathbb P_{(1,1,1,2)}$.
\end{lem}
\begin{proof}
First assume that $X$ is singular and let $\tilde X\rightarrow X$ be
an equivariant resolution of the singularities. By
\cite[cor.~3.6]{Mori82}, there exists a relative contraction
$\psi:\tilde X\rightarrow X'$ over $X$. Note that $\psi$ must be
divisorial. If $E$ is the exceptional divisor, use the classification
of \cite[thm.~3.3]{Mori82} to see that $E\cong \Ptwo$,
$\Pone\times\Pone$ or a singular quadric. As the map $S\rightarrow
\Aut(E)$ may not have a positive dimensional kernel, $S$ acts
transitively on $E$. This already rules out the singular quadric. No
$G$-invariant curve or divisor may intersect $E$. Consequently, the
$G$-exceptional set in $X'$ contains an isolated fixed point. By
\cite[thm.~1 on p.~113]{H-Oe}, $X'$ is a cone over a rational
homogeneous surface. Again using the \cite[thm.~3.3]{Mori82}, only
$\Pthree$ and the blow down of the $\infty$-section of $\mathbb
P(\mathcal O_{\Ptwo}(2)\oplus\mathcal O_{\Ptwo})$ are possible. This
variety is isomorphic to $\mathbb P_{(1,1,1,2)}$. By equality of the
\name{Picard}-numbers, $X'=X$.

If $X$ is smooth and homogeneous, then claim that $X\cong \Pthree$ or
$\QZ_3$. If the complement of the open orbit has dimension $<2$, then
use \cite[thm.~1 on p.~113 and thm.~1 on p.~121]{H-Oe} to yield the
claim (the other models occurring in the classification are either not
rational or have higher \name{Picard}-numbers). If the $G$-exceptional
set contains a divisor $E$, then $S$ acts non-trivially on $E$, and
$E\cong \Ptwo$ or $\Pone\times\Pone$. Now \cite[thms.~1 and 5]{Bad82}
apply, showing the claim.
\end{proof}

As a next step we rule the possibility out that the generic $S$-orbit
is a curve.
\begin{lem}
Under the assumptions of theorem~\ref{mainthm}, let $Y$ be a point. If
$G$ contains a subgroup $S\cong SL_2$, then the generic $S$-orbit is
of dimension 2 or 3.
\end{lem}
\begin{proof}
Assume to the contrary and let $C\subset X$ be any curve which is not
$SL_2$-invariant. Then $D:=\overline{SL_2.C}$ is an $S$-invariant divisor
and the generic $S$-invariant curve does not intersect $D$. A
contradiction to $X$ being minimal over a point!
\end{proof}

The case that the generic $SL_2$-orbit is 2-dimensional has been
classified in \cite{K98c}. The main result of that paper is that
$X\cong \QZ_3$,  $\Pthree$, $\mathbb P_{(1,1,1,2)}$ or $\mathbb
P_{(1,1,2,3)}$.

The case that $SL_2$ acts almost transitively will be considered
now. As a first step we recall that the assumption on $\QZ$-factorial
singularities already implies that $X$ is smooth.

\begin{lem}\label{factlem}
Let $X$ be a projective 3-dimensional variety with at most terminal
singularities, quasihomogeneous with respect to an algebraic action of
$SL_2$. Then either $X$ is smooth or not $\QZ$-factorial.
\end{lem}
\begin{proof}
If $X$ is not smooth, the singularities are isolated, hence fixed. Let
$p\in X$ be a singular point. Recall that $X$ can equivariantly
embedded into a projective space. Together with the complete
reducibility of $SL_2$ representations this yields an embedding of a
neighborhood $A$ of $p$ into a representation space $V$ such that $A$
is realized as the closure of an $SL_2$ orbit. A linearization
argument using the assumption that $SL_2$ has a three-dimensional
orbit yields that $p$ is necessarily the unique fixed point in $A$;
consequently, it's image is 0.

Follow the proof \cite[Lemma 5 on p.~210]{K85} in order to
construct two divisors $D_1$, $D_2 \subset X$ with $D_1 \cap  D_2 =
\{p \}$.
\end{proof}

Since $X$ is smooth, it must be contained in \name{Iskovskih}'s
list. It remains to find out which of the varieties in the list
actually occur in our context.

\begin{prop} 
In the setting of theorem~\ref{mainthm}, let $Y$ be a point. If
$S=SL_2$ acts almost transitively on $X$, then $X$ is one of the
following \name{Fano} varieties: $\Pthree$, $\QZ_3$, $V^S_{22}$ or
$V_5$.
\end{prop}

\begin{proof} We have already seen that $X$ is
necessarily smooth. Let $T<S$ be a maximal torus, and let $F\subset X$
be the set of $T$-fixed points.

If $F$ is discrete, use a linearization argument to see that the
\name{Lefschetz}-index of every fixed point is positive. The
\name{Borel} fixed point theorem asserts that $F$ is not empty. Thus,
$\chi (X)>0$ by the \name{Hopf} index theorem. We know already that
$b_0=b_6=1$, $b_1=b_5=0$ as $X$ is rational and $b_2=b_4=1$ by the
assumption that $\rho(X)=1$. Accordingly, $\chi(X)>0$ is possible iff
$b_3<4$. The classification of \name{Iskovskih} implies already that
only $\Pthree$, $\QZ_3$, $V_5$ and $V_{22}$ are possible. Recall that
the only quasihomogeneous variety of type $V_{22}$ is the special
$V^S_{22}$.

If $F$ is not discrete, then let $H$ be a component of $E$, the
complement of the open $S$-orbit in $X$, such that $\dim (F\cap
H)=1$. Since $S$ is 3-dimensional, $E$ is of pure dimension 2 and is
the support of an effective divisor generating the anticanonical
bundle $-K_X$. The $S$-action on $H$ cannot be almost transitive;
instead, the generic $S$-orbit must be 1-dimensional. Furthermore, $X$
does not contain an $S$-fixed point, or else a linearization at this
point would reveal a contradiction to the quasihomogeneous action of
$S$, there being no 3-dimensional representation of $SL_2$ with
3-dimensional orbits. This implies already that the normalization
$\tilde H$ of $H$ must be smooth. The closed and disjoint $S$-orbits
realize $\tilde H$ as a product $\tilde H \cong C\times\Pone$, where
$C$ is a smooth curve and $S$ acts on the second factor only. In
particular, there is no isolated $T$-fixed point in $\tilde H$, and
also none in $H$. As a next step, show that $H$ is smooth. In order
verify this claim, note that $F \cap H$ is smooth and that every
$S$-orbit in $H$ intersects $F \cap H$ transversally. If $H$ was
singular, let $x\in H_{sing}$ be a $T$-fixed point. If $U<SL_2$ is a
one-parameter group {\em not} fixing $x$, then by what we said above,
the map
\begin{eqnarray*}
F\times U &\rightarrow E \\ 
(f,u) &\mapsto u.f 
\end{eqnarray*}

has maximal rank at $(x,1)$, a contradiction! The adjunction formula
and the non-triviality of $K_H$ show that $H$ is \name{Fano}. So
$H\cong \Pone\times\Pone$ and \cite[thm.~5]{Bad82} yields the claim.
\end{proof}

\end{document}